# THIRTY-SIX UNSOLVED PROBLEMS IN NUMBER THEORY

by Florentin Smarandache, Ph. D.

University of New Mexico

Gallup, NM 87301, USA

**Abstract.**
Partially or totally unsolved questions in number theory and geometry especially, such as coloration problems, elementary geometric conjectures, partitions, generalized periods of a number, length of a generalized period, arithmetic and geometric progressions are exposed.

**Introduction.**

Mathematical philosophy?

The development of mathematics continues in a rapid rhythm, some unsolved problems are elucidated and simultaneously new open problems to be solved appear.

1. "Man is the measure of all things". Considering that mankind will last to infinite, is there a terminus point where this competition of development will end? And, if not, how far can science develop: even so to the infinite? That is . . .

The answer, of course, can be negative, not being an end of development, but a period of stagnation or of small regression. And, if this end of development existed, would it be a (self) destruction? Do we wear the terms of self-destruction in ourselves? (Does everything have an end, even the infinite? Of course, extremes meet.)

I, with my intuitive mind, cannot imagine what this infinite space means (without a beginning, without an end), and its infinity I can explain to myself only by means of a special property of space, a kind of a curved line which obliges me to always come across the same point, something like Moebus Band, or Klein Bottle, which can be put up/down (!)

I am not a specialist in physics, astronomy or philosophy, and I institute the infinite only theoretically --from the point of view of mathematics (of course, unilaterally).

2. Mathematics knows a high degree of abstraction, idealization, and generalization. And I ask, is it possible to build a pure mathematical model for society? You will answer, "it would be too rigid". And you are right because the non-elastic systems stop the progress. But I should replay, "they would be based upon logic". In the future could we become human robots, having the time programmed to the second (or even thousandth of a second!), elements of a mathematical-cybernetic system?

3. How do you imagine mathematics over 1,000 years? What about 1,000,000 year? (I exaggerate, of course.) What

other new branches will appear? (Some will be ended, out of date?) (I'm not a futurist [Toepler]).

You may consider these questions too general, or too easy, but what can be more difficult than trying to say a maximum of ideas with a minimum of words? You are right, sometimes a too general idea becomes a common one. Maybe you think that asking questions is easy, but let me contradict you. (R. K. Guy said that asking questions is an art.) And after all, aren't the theories born, at their turn, from questions? (Maybe in this essay the questions are too direct, you are right again.)

4. If we consider "Math (t)", the development of mathematics at a time "t" (considered from the appearance of life on Earth) then

$$(\exists) \; L = \lim_{t \to \infty} \text{Math}(t) \; ?$$

And if it is yes, what is it equal to? To $\infty$? In case of total self-destruction should we obtain L = 0? And if life would reappear and the development would start again then should we draw the conclusion that $(\nexists) \; L$? (cyclical development).

5. In case of a total (self) destruction and reappearance of life on Earth, how would they call the Pythagoras' theorem, Euclid's geometry, Birkoff's axioms,

Erdös's open problems, and so on?  Will mankind pass through the same phases of development?  Or, if it would exist, another civilization at the same time, how should they call these results?  And how should we call them if these two (or more) civilizations united?  (I have arrived in the field of history and of mathematical philosophy, which is not the purpose of this paper.)  (All these questions can be also extended to other fields of knowledge.)

I can imagine computers with fantastic memories having the whole mathematics divided like a mosaic:  this theorem belongs to X, the other one to Y, this sentence belongs to everybody, the other sentence belongs to nobody--the one who will invent has not been born yet, but he will be born!  A real dictionary of names and ideas, science divided in a finite (but, however, infinite) number of cells, each of them having a strict delimitation, with its own history, and the future generations will add new and new cells.

Although the applied mathematics, the integral calculus, the operator theory are the queens, the primitive arithmetic still fascinates the world because of its apparent elementary problems--very easy to be understood by everybody, but … .

Why is mankind still interested in these easy problems, of mathematical juggler?  I think that it is possible thanks to their simplicity in exposure.  ("Ah, it's easy", one says, and to solve it you discover that you plunge into a

labyrinth. And, hence, appears a paper: "On a conjecture . . .", "On the unsolved problem . . ." etc.)

I am sure that the "unsolved" problems presented in these pages will be (or have already been before the appearance of this essay) easy for many mathematicians, but for me they were an obsessions. W. Sierpiński was optimistic when he said that if mankind lasted then all these unsolved problems would be solved.

All the material in this paper is original in the author's opinion. He wanted to gather in this paper a variety of material, for the sake of an harmony of contraries.

"I want to be a mathematician", as P. R. Halmos, and for the I began to play: rebus + mathematics, literature + mathematics, and even rebus + literature! So, please, do not wonder at this essay.

Springer, 1984.

## UNSOLVED PROBLEM: 1

Find all integer sequences $\{a_n\}_{n \in \mathbb{N}^*}$ defined as follows:

(I) $(\forall) i \in \mathbb{N}^*$, $(\exists) j, k \in \mathbb{N}^*$, $i \neq j \neq k \neq i$, such that

$a_i \equiv a_j \pmod{a_k}$.

(II) $(\forall) i \in \mathbb{N}^*$, $(\exists) j, k \in \mathbb{N}^*$, $i \neq j \neq k \neq i$, such that

$a_j \equiv a_k \pmod{a_i}$.

## UNSOLVED PROBLEM: 2

Let $d > 0$. Questions:

(a) What is the maximum number of points included in a plane figure (generally: in a space body) such that the distance between any two points is greater than or equal to d?

(b) What is the minimum number of points $\{A_1, A_2, \ldots\}$ included in a plane figure (generally: in a space body) such that if it includes another point A then there would be an $A_i$ with $AA_i < d$?

## UNSOLVED PROBLEM: 3

(a) Let $a_1, \ldots, a_n$ be distinct digits of the set $\{0, 1, 2, \ldots, 9\}$, for a given n, $1 \leq n \leq 9$. How many distinct primes can we make up with all these digits? More generally: when $n \in \mathbb{N}^*$ and $a_1, \ldots, a_n$ are distinct positive integers.

(b) Let $a \in \{0, 1, \ldots, 9\}$. How many digits of a does the n-th prime contain? But n! ? But $n^n$ ? More generally: when $a \in \mathbb{N}$.

**Comment**

"The sizes $P_n$, n!, $n^n$ have jumps when $n \to n + 1$, hence the analytical expressions are approximate only. Moreover, the results depend on the <u>exact</u> (and not approximate) value of these sizes" (E. Grosswald [1]).

"(a) can be solved quickly on a modern computer" (R. K. Guy [2]).

# UNSOLVED PROBLEM: 4

Rationalize the following fraction:

$$1 \bigg/ \sum_{i=1}^{n} a_i^{\wedge}(1/k_i).$$

# UNSOLVED PROBLEM: 5

Mathematical Logic:

Is it true that for any question there is at least an answer? Reciprocally: Is any assertion the result of at least a question?

# UNSOLVED PROBLEM: 6

Is it possible to construct a function which obtains all irrational numbers? But all transcendental numbers?

**UNSOLVED PROBLEM: 7**

Given n points in space, four by four non-coplanar, find a maximum m having the property that there are m points among n ones that constitute the vertexes of a convex polyhedron. [An extension of the following conjecture: No matter how one chooses $2^{m-2} + 1$ points in plane, three by three non-collinear, there are among these m joints which are the vertexes of a convex polygon. (Ioan Tomescu, Problems of combinatory and graph theory [Romanian], Bucharest, EDP, 1983.) For m = 5 the conjecture was proved; it was still proved that it can choose $2^{m-2}$ points in plane, three by three non-collinear, such that any m ones among these do not constitute the vertexes of a convex polygon.]

**UNSOLVED PROBLEM: 8**

What is the maximum number of circles of radius 1, at most tangential by twos, which are included into a circle of radius n? (Gamma 1/1986). This problem was generalized by Mihaly Bencze, who asks the maximum number of circles of radius $\emptyset(n)$, at the most tangential by twos, which are included into a circle of radium n, where $\varphi$ is a function of n (Gamma 3/1986).

Study a similar problem for circles of radius 1 included into a given triangle (on Malfatti's problem). Similar questions for spheres, cones, cylinders, regular pyramids, etc. More generally: planar figures included into a given planar figure. And in the space, too.

## UNSOLVED PROBLEM: 9

(a) Let $m \geq 5$ an integer. Find a minimum n (of course, n depends on m) having the property: no matter how one chooses n points in space, four by four non-coplanar, there exist m ones among these which belong to a surface of a sphere.

(b) Same question for an arbitrary spatial figure (for example: cone, cube, etc.).

(c) Similar problems in plane (for $m \geq 4$, and the points: three by three non-collinear).

## UNSOLVED PROBLEM: 10

Let $a_1, a_2, \ldots, a_m$ be digits. Are there primes, on a base b, which contain the group of digits $\overline{a_1 \ldots a_m}$ into its writing? (For example, if $a_1 = 0$ and $a_2 = 9$ there are primes as 1$\underline{09}$, 4$\underline{09}$, 7$\underline{09}$, 8$\underline{09}$, ...) But n! ? But $n^n$ ?

## UNSOLVED PROBLEM: 11

**Conjecture**

Let $k \geq 2$ a positive integer. The diophantine equation:

$$y = 2 x_1 x_2 \ldots x_k + 1$$

has an infinity of solutions of primes. (For example: $571 = 2 \cdot 3 \cdot 5 \cdot 19 + 1$, $691 = 2 \cdot 3 \cdot 5 \cdot 23 + 1$, or $647 = 2 \cdot 17 \cdot 19 + 1$, when k = 4, respectively, 3). (Gamma 2/1986).

## UNSOLVED PROBLEM: 12

Let $d_n$ be the distance between two consecutive primes, $d_n = \frac{1}{2}(p_{n+1} - p_n)$, n = 1, 2, ... . Does $d_n$ contain an infinite number of primes?

Does $d_n$ contain numbers of the form n! ? But of the form $n^n$ ?

Let $i \in \mathbb{N}^*$, and $d_n^{(i)} = \frac{1}{2}(p_{n+i} - p_n)$, and $d_{n,i}^{(i-1)} =$
$= \frac{1}{2}(d_{n+1} - d_n)$ for $n = 1, 2, \ldots$

The same questions.

(Gamma 2/1986).

**UNSOLVED PROBLEM:  13**

**Conjecture:**

No matter how the points of a plane are colored with n colors, there exists a color which fulfills all distances [i.e., $(\exists)$ a color c, $(\forall)$ d $\geq$ o, $(\exists)$ the points A, B colored in c, such that the line AB is colored in c and $AB = d$].

(The result would implicitly be true in space, too.)

**UNSOLVED PROBLEM:  14**

Let $k, n \in \mathbb{N}^*$, $k < n$. We define a decreasing sequence of integers: $n_0 = n$ and $n_{i+1} = \max \{p, p \, n_i - k, p \text{ is a prime}\}$, for $i \geq 0$.

Find the length (the number of terms) of this sequence (Gamma 2-3/1987).

## UNSOLVED PROBLEM: 15

Spreading to infinite van der Waerden's theorem: Is it possible to partition $N^*$ into an infinity of arbitrary classes such that at least one class contain an arithmetic progression of $\ell$ terms ($\ell \geq 3$)?

Find a maximum $\ell$ having this property.

## UNSOLVED PROBLEM: 16

Let $a \in Q \setminus \{-1, 0, 1\}$. Solve the equation:

$$xa^{\frac{1}{x}} + \frac{1}{x}a^x = 2a.$$

[A generalization of the problem 0:123, Gazeta Matematicā, No. 3/1980, p. 125.]

# UNSOLVED PROBLEM: 17

(a) If $(a, b) = 1$, how many primes does the progression $ap_n + b$, $n = 1, 2, \ldots$, contain? where $p_n$ is the n-th prime. But numbers of the form $n!$ ? But $n^n$ ?

(b) Same questions for $a^n + b$, $a \notin \{\pm 1, 0\}$.

(c) Same questions for $k^k + 1$ and $k^k - 1$, $k \in \mathbb{N}^*$.

(Gamma 2/1986)

# UNSOLVED PROBLEM: 18

(a) Let n be a non-null positive integer and $d(n)$ the number of positive divisors of n. Of course, $d(n) \leq n$, and $d(n) = 1$ if and only if $n = 1$. For $n \geq 2$ we have $d(n) \geq 2$. Find the smallest k such that

$$\underbrace{d(d(\ldots d(n)\ldots))}_{k \text{ times}} = 2$$

(b) Let $\sigma(n) = \sum_{\substack{d/n \\ d>0}} d$ and m a given positive integer

Find the smallest k such that

$$\underbrace{\sigma(\sigma(\ldots \sigma(2)\ldots))}_{k \text{ times}} \geq m.$$

# UNSOLVED PROBLEM: 19

Let $a_1, a_2, \ldots$ be a strictly increasing sequence of positive integers, and $N(n)$ the number of terms of the sequence not greater than n.

(1) Find the smallest k such that

$$\underbrace{N(N(\ldots N(n)\ldots))}_{k \text{ times}} \text{ is constant, for a given } n.$$

(2) If $a_1 \notin \{0,1\}$, find the smallest k such that

$$a_k \geq n, \text{ for a given } n.$$

```
        a
    t      .
     i       .
      m       .
       e       a
        s       1
```

Particular Cases:

    (a)   When $\{a_i\}_i$ is the sequence of primes.

    (b)   When $\{a_i\}_i$ is the sequence of m-th powers, for a given m : $0, 1, 2^m, 3^m, \ldots$ .

    For example, the sequence of perfect squares.

    (c)   $\{a_i\}_i$ is anyone of the well known sequences of positive integers, for example Fibonacci (or Lucas, Fermat, Bernoulli, Mersenne, Smarandache, van der Waerden, etc.) numbers.

## UNSOLVED PROBLEM: 20

    Let k be a non-zero integer. There are only a finite number of solutions in integers p, q, x, y, each greater than 1, of the equation $x^p - y^q = k$.

    (On Catalan's conjecture) [For k = 1 this was conjectured by Cassels (1953) and proved by Tijdeman (1976).]

(Gamma 2/1986)

## UNSOLVED PROBLEM: 21

Let $\{x_n\}_{n \geq 1}$ be a sequence of integers, and $0 \leq k \leq q$ a digit. We define a <u>sequence of position</u>:

$$U_n^{(k)} = U^{(k)}(X_n) = \begin{cases} i, & \text{if } k \text{ is the i-th digit of } x_n; \\ 0, & \text{in the other ones} \end{cases}$$

(For example: if $x_1 = 5$, $x_2 = 17$, $x_3 = 715$, ... and $K = 7$, then

$$U_1^{(7)} = U^{(7)}(x_1) = 0, \quad U_2^{(7)} = 2, \quad U_3^{(7)} = 1, \ldots)$$

(1) Study $\{U^{(k)}(p_n)\}_n$, where $\{p_n\}_n$ is the sequence of primes. Convergence, monotony.

The same questions for the sequences:

(2) $x_n = n!$, $n = 1, 2, \ldots$

(3) $x_n = n^n$, $n = 1, 2, \ldots$

Do the sequences of position $U_n^{(k)}$ contain an infinity of primes?

More generally: when $\{x_n\}_n$ is a sequence of rational numbers and $k \in \mathbb{N}$.

**UNSOLVED PROBLEM: 22**

Let M be a number in a base b. All distinct digits of M are named <u>generalized period</u> of M. (For example, if M = = 104001144, its generalized period is $g(M) = \{0, 1, 4\}$.). Of course, $g(M) \subseteq \{0, 1, 2, ..., b - 1\}$.

The <u>number of generalized period</u> of M is equal to the number of groups of M such that each group contain all digits of $g(M)$. (For example, $n_g(M) = 2$, M = $\underline{104}\ \underline{001144}$).
$\phantom{xx}1\phantom{xx}2$

<u>Length of the generalized period</u> is equal to the number of its digits. (For example, $\ell_g(M) = 3$.)

**Questions:**

(1) Find $n_g$, $\ell_g$ for $p_n$, $n!$, $n^n$, $\sqrt[n]{n}$.

(2) For a given $k \geq 1$, is there an infinity of primes $p_n$, or $n!$, or $n^n$, or $\sqrt[n]{n}$ which have a generalized period of length k? But which have the number of generalized periods equal to k?

(3) Let $a_1, a_2, ..., a_h$ distinct digits. Is there an infinity of primes $p_n$, or $n!$, or $n^n$, or $\sqrt[n]{n}$ which have as generalized period the set $\{a_1, a_2, ..., a_h\}$?

**Remark:**

There exist arithmetic and geometrical progressions which contain an infinity of terms of given generalized period. For example, if g = {6, 7}, we construct an arithmetic progression: <u>67</u>, 1<u>67</u>, 2<u>67</u>, ... and a geometrical one: <u>67</u>, <u>67</u>0, <u>67</u>00, ...

## UNSOLVED PROBLEM: 23

Find the maximum r such that: the set $\{1, 2, \ldots, r\}$ can be partitioned into n classes such that no class contains integers x, y, z with $xy = z$.

(On Schur's Problem)

Same question when $x^y = z$.

Same question when no integer can be the sum of another integer of its class. (A generalization of Schur's Problem.)

## UNSOLVED PROBLEM: 24

Let $N = \{1, 2, \ldots, n\}$. Find the maxim number of elements extracted from N such that any m from these be not an arithmetic progression ($n > m > 2$).

Same question when the m elements must not be a geometrical progression.

More generally: Let R be a given m - ary relation on a field N. What is the maximum number of elements extracted from N such that any m from these be not in this relation? What happens when N has continuation power (for example when N is an interval from R)? [On van der Waerden's theorem.]

**UNSOLVED PROBLEM: 25**

Let $ù$ numbers be $a + bù$, where $ù$ is a complex n-th root of unity, $ù^{n-1} + ù^{n-2} + \ldots + 1 = 0$, which enjoy unique factorization. The units are: $\pm 1, \pm ù, \pm ù^2, \ldots, \pm ù^{n-1}$.

**Conjecture:**

The configuration of $ù$ primes are symmetric of the 2n regular polygon.

(On Gaussian primes. A generalization of Einstein's integers.)

**Reference:**

[1] R. K. Guy, Unsolved Problems in Number Theory, Springer-Verlag, New York, Heidelberg, Berlin, 1981, pp. 20-22, A16.

**UNSOLVED PROBLEM: 26**

The equation $x^3 + y^3 + z^3 = 1$ has as solutions (9, 10, -12) and (-6, -8, 9). How many other nontrivial integer solutions are there?

**References:**

[1] V. L. Gardiner, R. B. Lazarus, and P. R. Stein, Solution of the Diophantine Equation $x^3 + y^3 = z^3 - d$,

Math. Comput. 18 (1964) 408-413; MR 31 #119.

[2] J. C. P. Miller and M. F. C. Woollett, Solutions of the Diophantine Equation $x^3 + y^3 + z^3 = k$, J. London: Math. Soc. 30 (1955) 101-110; MR 16, 979.

[3] R. K. Guy. Unsolved Problems in Number Theory. Springer-Verlag, New York, Heidelberg, Berlin, 1981, p. 84, D5.

## UNSOLVED PROBLEM: 27

Daniel Silverman asked if $\prod_{n=1}^{m} \frac{p_n+1}{p_n-1}$, where $p_n$ is the n-th prime, is an integer for others $m \in \{1, 2, 3, 4, 8\}$. We conjecture that

$$R_m = \prod_{n=1}^{m} \frac{p_n+k}{p_n-k}, \text{ with } K \in \mathbb{N}^*,$$

is an integer for a finite number of values of m. There is an infinite number of k for which no $R_m$ is an integer.

**Reference:**

[1] R. K. Guy, Unsolved Problems in Number Theory, Springer-Verlag, New York, Heidelberg, Berlin, 1981, p. 57, B48.

# UNSOLVED PROBLEM: 28

## (ON A PROBLEM WITH INFINITE SEQUENCES)

Let $1 \leq a_1 < a_2 < \ldots$ be an infinite sequence of integers such that any three members do not constitute an arithmetical progression. Is it true that always $\sum_{n \geq 1} 1/a_n \leq 2$? Is the function

$$S(\{a_n\}_{n \geq 1}) = \sum_{n \geq 1} 1/a_n$$

bijective?

For example, $a_n = p^{n-1}$, $n \geq 1$, $p$ is an integer $> 1$, has the property of the assumption, and $\sum_{n \geq 1} 1/a_n = 1 + \frac{1}{n-1} \leq 2$. Analogously for geometrical progressions.

More generally: let $f$ be a function $f: R_+^m \to R_+^*$. We construct a sequence $0 < a_1 < a_2 < \ldots$ such that there is no $(a_{i_1}, \ldots, a_{i_m}, a_{i_{m+1}})$ with $f(a_{i_1}, \ldots, a_{i_m}) = a_{i_{m+1}}$. Find

$$\max_{\{a_n\}_{n \geq 1}} \sum_{n \geq 1} 1/a_n .$$

(It's a generalization of a question from the problem E28, R. K. Guy, Unsolved Problems in Number Theory, Springer-Verlag, 1987, p. 127.)

Is the function

$$S(\{a_n\} n \geq 1)_{n \geq 1} = \sum 1/a_n$$

bijective?

# UNSOLVED PROBLEM: 29

Let ó(n) be sum of divisors of n, $\prod(x)$ number of primes not exceeding x, ù(n) number of distinct prime factors of n, d(n) number of positive divisors of n, p(n) the largest prime factor of n.

Let $f^{(k)}$ note $\underbrace{f \circ f \circ \ldots \circ f}_{k \text{ times}}$, for all function f.

Find the smallest k for which:

(a) for fixed n and m, we have $ó^{(k)}(n) > m$.

(b) for a fixed real x, with $x \geq 2$, we have $\prod^{(k)}(x) = 1$.

(c) for a fixed n, we have $ù^{(k)}(n) = 1$.

(d) for fixed n and m, we have $d^{(k)}(n) > m$.

(e) for a fixed n, we have

$$\underbrace{p(p(\ldots(p(n)-1)\ldots)-1)-1}_{k \text{ times}} = 1.$$

Generalize for $ó_r(n)$, $\prod(x; a, b)$, $Ù(n)$.

# UNSOLVED PROBLEM: 30

# AN EXTENSION OF CARMICHAEL'S CONJECTURE

**Conjecture:**

$\forall a \in \mathbb{N}^*, \forall b \in \mathbb{N}^*, \exists c \in \mathbb{N}^*: \varphi(a) \cdot \varphi(b) = \varphi(c)$. For $a = 1$ it results Carmichael's conjecture. If this conjecture is true, then by mathematical induction it finds:

$\forall a_1, \ldots, a_n \in \mathbb{N}^*, \exists b \in \mathbb{N}^*: \varphi(a_1) \ldots \varphi(a_n) = \varphi(b)$.

**Reference**

[1] R. K. Guy, Unsolved Problems in Number Theory, Springer-Verlag, 1981, p. 53, B39.

# UNSOLVED PROBLEM: 31

# (ON CRITTENDEN AND VANDEN EYNDEN'S CONJECTURE)

Is it possible to cover all (positive) integers with n geometrical progressions of integers?

*Find a necessary and sufficient condition for a general class of positive integer sequences such that, for a fixed n, there are n (distinct) sequences of this class which cover all integers.

**Comment:**

(a) No. Let $a_1, \ldots, a_n$ be respectively the first term of each geometrical progression, and $q_1, \ldots, q_n$

respectively their rations.  Let ñ be a prime, ñ ∉ {$a_1$, ..., $a_n$, $q_1$, ..., $q_n$}.  Then ñ does not belong to the union of these n geometrical progressions.

(b)*  For example, the class of progressions:

$$A_f = \{\{a_n\}_{n \geq 1} : a_n = f(a_{n-1}, \ldots, a_{n-i}) \text{ for } n \geq i + 1, \text{ and } i, a_1, \ldots, a_i, \in \mathbb{N}^*\}$$

with the property:

$$\exists y \in \mathbb{N}^*, \forall (x_1, \ldots, x_i) \in \mathbb{N}^{*i} : f(x_1, \ldots, x_i) \neq y,$$

does it cover all integers?

But, if we change the property:

$$\forall y \in \mathbb{N}^*, \exists (x_1, \ldots, x_i) \in \mathbb{N}^{*i} : f(x_1, \ldots, x_i) = y?$$

(Generally no; see the geometrical progressions.)

This (solved and unsolved) problem remembers Crittenden and Vanden Eynden's conjecture.

### UNSOLVED PROBLEM: 32

Consider the following equation:

$$(a - b \sqrt[m]{n}) x + c \sqrt[m]{n} \cdot y + \sqrt[p]{q} \cdot z + (d + e \cdot w) \sqrt[r]{s} = 0,$$

where a, b, c, d, e are constant integers; and the m-th, p-th and r-th roots are irrational distinct numbers. What conditions must the parameters m, n, p, q, r and s accomplish such that the equation admits integer solutions (x, y, z and w being variables)?

**UNSOLVED PROBLEM: 33**

Find the maximum number of interior points inside of a parallelogram having an angle of ð/3 such that the distance between any two points is greater than or equal to 1. (The same question for a prism where all the faces are parallelograms with an angle of ð/3.)

More generally: let d > o. Questions:

(a) Which is the maximum number of points included in a plane figure (a space body) such that the distance between any two points is greater or equal than d?

(b) Which is the minimum number of points $\{A_1, A_2, ...\}$ included in a plane figure (a space body) such that if another point A is added then there is an $A_i$ with $AA_i$ < d?

(c) Other variants of these questions if it considers:

 (1) only interior points; (2) interior and frontier points; (3) points on frontier only; (4) the distance is strictly greater than d; and (5) $AA_i$ ≤ d.

**Solution:**

(I) We consider an equilateral triangle's network as indicated in the diagram displayed below, where $\ell_1$ and $\ell_2$ are the sides of our parallelogram. Clearly, this network gives the optimum construction of interior and frontier points keeping our property. It finds $(\lfloor \ell_1 \rfloor + 1)(\lfloor \ell_2 \rfloor + 1)$ points, where $\lfloor x \rfloor$ is the greatest integer less or equal than x. If $\ell_1, \ell_2 \in N^*$ we cannot take more than $\ell_1$ interior points on a rule and $\ell_2$ interior points on the other one, because it is not permitted to take points on the frontier.

If, for example, $\ell_1 \in N$, we can take $\lfloor \ell_1 \rfloor + 1$ interior points on a side. In conclusion: $n_{max} = \lceil \ell_1 \rceil \bullet \lceil \ell_2 \rceil$ interior points, where $\lceil x \rceil$ is the little integer greater or equal than x.

(II) For the prism of our problem, having the sides $\ell_1, \ell_2, \ell_3$, of course, $n_{max} = \lceil \ell_1 \rceil \bullet \lceil \ell_2 \rceil \bullet \lceil \ell_3 \rceil$ [it results from (I) by considering the parallelograms $(\ell_1, \ell_2)$ and $(\ell_2, \ell_3)$].

(III) These are generally open (unsolved) questions. For particular cases see [1].

**Reference:**

[1] Smarandache, Florentin, Problèmes avec et sans ... problèmes! (problems 5.43 [p. 67], respectively 5.44 [p. 62]), Somipress, Fès, Morocco, 1983 (M. R.: 84R: 00003).

*Find all real solutions of the equation $x^y - \lfloor x \rfloor = y$, where $\lfloor x \rfloor$ is the greatest integer less than or equal to x.

Solution: It is a generalization of a Putnam competition problem: ($x^3 - \lfloor x \rfloor = 3$).

(1) When $y \in \mathbb{R} \setminus \mathbb{Q}$ the author is not able to answer.

(2) If $y = 0$ then $x \in [1, 2)$.

(3) If $y = 1$, $\nexists x \in \mathbb{R}$.

(4) If y is an odd integer > 1, then $x = \sqrt[y]{y+1} \in (1, 2)$, is the unique solution of our equation.

Let's use the functions:

$f(x) = x^n - n$, $n \in \mathbb{N}$, $n > 1$, where n is an odd integer,

and $g(x) = \lfloor x \rfloor$, $f, g: \mathbb{R} \to \mathbb{R}$.

(5) If $y = 2$, then $x_1 = \sqrt{3}$ and $x_2 = -1$.

(6) If $y \in \mathbb{N}$ and y is an even integer $\geq 4$, then $x_1 = \sqrt[y]{y+1}$, $x_2 = -\sqrt[y]{y-2}$. (It is sufficient to observe that $f(x) = x^n - n$, where n is an even integer $\geq 4$, is an even function; we consider $f: \mathbb{R} \to \mathbb{R}$ (Fig. 2).).

(7) If $y = -1$ the $x = -1/2$, because it results in $1/x \in \mathbb{Z}^*$, hence $x = 1/k$ with $k \in \mathbb{Z}^*$, whence $k = -2$.

(8) If $y \in \mathbb{Z}$, $y < -1$, $\nexists x \in \mathbb{R}$, because it results $1/x^m \in \mathbb{Z}^*$, whence $x^m = 1/k$, $k \in \mathbb{Z}^*$, where $m = -y > 0$. For $\sqrt[m]{\phantom{x}}$ $\sqrt[m]{\phantom{x}}$ an even $m$, we have $x = \pm \sqrt{1/k}$ with $k \geq 1$, but $\lfloor \pm \sqrt{1/k} \rfloor \in \{0, \pm 1\}$ and $k + m \geq 2$; hence $\nexists\ x \in R$. For an odd $m$, we have $x = \sqrt[m]{1/k}$ with $k \in Z^*$; hence $k - \lfloor \sqrt[m]{1/k} \rfloor = -m$ whence $k = -m - 1 < 0$ and $x = -\sqrt[m]{1/(m+1)} = -(1-y)^{1/y}$.

(9) Let $y = 1/n$, $n \in N$, $n$ is odd $\geq 3$, then $x_1 = (1/n)^n$, $x_2 = (\overline{n}^{n-1})^n$. The function $f(x) = x^{\frac{1}{n}} - \frac{1}{n}$ $\in N$, $n$ is odd $\geq 3$, $f: R \to R$.

For $x \geq 2$, $h(x) = x^{\frac{1}{n}} - x + 1 - \frac{1}{n} < 0$ (because $h'(x) < 0$ when $x \geq 2$, hence $h(x) \leq h(2) < 0$ for $x \geq 2$). There exists an unique positive real solution.

For $x \leq -1$, $k(x) = x^{\frac{1}{n}} - \frac{1}{n} - x > 0$; there exists an unique negative real solution.

(10) If $y = 1/n$, $n \in N$, $n$ is even $\geq 2$, then $x \geq 0$ involves $x = (1/n)^n$, because $f$ is an even function (Fig. 4).

(11) If $y = -1/n$, $n \in N$, $n \geq 2$.

$$1/x^{1/n} + 1/n = \lfloor x \rfloor.$$

Whence $x^{1/n} \in Q$, it results $x = a^n$, $a \in Q$. Hence $1/a + 1/n \in Z$, whence $1/a = k - 1/n$, $k \in Z$, thus $a = \frac{n}{kn-1}$ or $k = \lfloor \frac{n}{kn-1} \rfloor$. Of course $k \neq 0$. If $k = 1$, we have $1 + \frac{1}{\left(kn-1\right)^n} > 2$. If $k \geq 2$, we have $\frac{n^n}{kn-1} < 1$.

(a) If n is odd then $k = -1$, whence $x = -\dbinom{n-1}{\frac{n+1}{2}}^n \bigg/ \dbinom{kn-1}{n}^n$ is a solution.

(b) If n is even, there is no $x \in R$. For $k \leq -2$, we have $\dfrac{n^n}{\binom{kn-1}{m}} > -1$.

(12) $y = \dfrac{\binom{kn-1}{m}}{n} \in Q\backslash Z$, $m \neq 1$.

* The author is not able to solve the equation in this case.

### UNSOLVED PROBLEM: 35

Prove that on a circular disk of radius r there are least n points such that the distance between any two is greater or equal than d, where:

$$n = \sum_{k \geq 0}^{\left[\frac{r}{d}\right]} [\eth/\arcsin d/(2(r-kd))] \text{ if } \frac{r}{d} \notin N,$$

or

$$n = 1 + \sum_{k \geq 0}^{\frac{r}{d}-1} [\eth/\arcsin d/(2(r-kd))] \text{ if } \frac{r}{d} \in N.$$

Is n equal to the maximum number of points on the disk with this property?*

Generalize for an arbitrary figure in plane.**

Generalize also for an arbitrary corps in space.**

**Proof:**

(a) Let $æ_d(r)$ be the maximum number of points on the circumference of a circle of radius r such that the distance between any two is greater or equal than d. Hence, the cord which unites two points from these is $\geq$ d (see the picture 1). We take it precisely d.

$$\frac{d}{\sin x} = \frac{r}{\sin(\frac{\pi}{2}-x)} \quad \text{whence} \quad \sin\frac{x}{2} = \frac{d}{2r} \quad \text{hence}$$

$$x = 2\arcsin\frac{d}{2r} \quad \text{(in radians). We divide } 2\pi \text{ to } x$$

and it results in:

$$æ_d(r) = [\pi/\arcsin\frac{d}{2r}].$$

We proceed analogously with the circle $C_1$ of radius r-d, concentric with the first $C_0$, obtaining $æ_d(r-d)$, etc. This method ends at the step $k = [\frac{r}{d}]$ for which $0 < r-kd < d$. When $\frac{r}{d}$ is an integer, the last drawn circle will be a point, precisely the center of these circles.

This construction mode (with point networks lying on concentric circles such that between two some circles the distance is equal to d) ensures the distance condition of all points.

It remarks that if since start d > 2r the arcsin there is not, therefore our problem is impossible (there is no point). And, if d = 0, we obtain an infinite of points.

(b) *This construction is close to an optimum one (in the author's conception). But the author cannot prove if this is or is not optimum. There are many constructions of point networks [on squares, on (equilateral) triangles, etc.]. For which point networks is n maximum? (Here it is an open question).

In our problem, when $\frac{r}{d}$ is very great perhaps the following <u>construction</u> is more advantageous (see picture 2):

We take a point $P_1$ on the circumference $F_0$ of our disk $Z_0$. With a compass we draw a circle arc $A_1$ (of radius d, having the center in $P_1$) which cuts $F_0$ in $P_2$ (to right). Afterwards we again draw a circle arc $A_2$ (of radius d, having the center in $P_2$) which cuts $F_0$ in $P_3$ (to right), etc. On $F_0$ we find $æ_0$ points.

We can still take another points in the shaded zone $Z_1$ only. We construct these points on the frontier $F_1$ of $Z_1$, analogous: with a compass of radius d, of center $R_1$ at the start, etc. ($R_1 = A_1 \cap A_\ell$, where $A_\ell$ is the last circle arc). On $F_1$ we find $æ_1$ points.

This method ends when $Z_k = \emptyset$. It obtains at least n points having the property of our assumption. But, does it obtain more points than n?

(c)  **These are two general open questions.  For particular cases see [1].

**UNSOLVED PROBLEM:   36**

(An Hypothesis Extending the ERP-Paradox)

a) Foreword.
What's new in science (physics)?
According to researchers from the University of Innsbruck in Austria (December 1997):
- photon is a bit of light, the quantum of electromagnetic radiation (quantum is the smallest amount of energy that a system can gain or lose);
- polarization refers to the direction and characteristics of the light wave vibration;
- if one uses the entanglement phenomenon, in order to transfer the polarization between two photons, then: whatever happens to one is the opposite of what happens to the other; hence, their polarizations are opposite of each other;
- in quantum mechanics, objects such as subatomic particles do not have specific, fixed characteristic at any given instant in time until they are measured;
- suppose a certain physical process produces a pair of entangled particles A and B (having opposite or complementary characteristics), which fly off into space in the opposite direction and, when they are billions of miles apart, one measures particle A; because B is the opposite , the act of measuring A instantaneously tells B what to be; therefore those instructions would somehow have to travel between A and B faster than the speed of light; hence, one can extend the Einstein-Podolsky-Rosen paradox and Bell's inequality and assert that the light speed is not a speed barrier in the universe.

b) Scientific Hypothesis:
We even promote the hypothesis that: there is no speed barrier in the universe, which would theoretically be proved by increasing, in the previous example, the distance between particles A and B as much as the universe allows it, and then measuring particle A.

c) An Open Question now:
If the space is infinite, is the maximum speed infinite?

Kishinev, December 1994.

[10] Smarandache, Florentin, "There Is No Speed Barrier In The Universe", <Bulletin of Pure and Applied Sciences>, Delhi, India, Vol. 17D (Physics), No. 1, p. 61, 1998.

[11] Smarandache, Florentin, "There Is No Speed Barrier In The Universe",
http://www.gallup.unm.edu/~smarandache/NoSpLim.htm,
http://www.geocities.com/m_l_perez/SmarandacheHypothesis/Sm-Hyp.htm (in Quantum Physics online journal).

[12] Suplee, Curt, " 'Beaming Up' No Longer Science Fiction", <Albuquerque Journal>, December 11, 1997.

13] Walorski, Paul (A.B. Physics), Answer to J. Gilbert, Ask Experts: http://www.physlink.com/ae86.cfm.